\documentclass[12pt]{article}
\usepackage{amsmath,amsfonts,amssymb,amsthm}
\usepackage{bbm}
\usepackage{graphicx}
\usepackage{stmaryrd}
\usepackage{enumerate}
\usepackage{mathtools}
\usepackage[font=footnotesize]{caption}

\usepackage[usenames,dvipsnames,svgnames]{xcolor}
\usepackage{relsize}

\newcommand{\edge}{\mathsmaller{\mid}}

\newcommand{\grad}{\nabla}

\newcommand{\laplace}{\Delta}

\def\Xint#1{\mathchoice
{\XXint\displaystyle\textstyle{#1}}%
{\XXint\textstyle\scriptstyle{#1}}%
{\XXint\scriptstyle\scriptscriptstyle{#1}}%
{\XXint\scriptscriptstyle\scriptscriptstyle{#1}}%
\!\int}
\def\XXint#1#2#3{{\setbox0=\hbox{$#1{#2#3}{\int}$ }
\vcenter{\hbox{$#2#3$ }}\kern-.6\wd0}}
\def\avint{\Xint-}
\def\aviint{\avint\!\!\!\avint}

\renewcommand{\div}{\grad\cdot}

\newcommand{\R}{\mathbf{R}}

\newcommand{\Ha}{\ensuremath{\mathcal{H}}}
\newcommand{\D}{\ensuremath{\mathcal{D}}}

\usepackage{mathtools}
\mathtoolsset{showonlyrefs,showmanualtags}

\newtheorem{theorem}{Theorem}

\newtheorem{example}{Example}

\newcommand{\noto}{\,\,\not\!\!\longrightarrow}

\newcommand{\tacka}{\, \cdot\,}

\usepackage{suffix}
\usepackage{mathtools}
\usepackage{hyperref}

\newcommand{\cT}{\mathcal{T}}

\begin{document}

\title{Optimal stability estimates for continuity equations}
\author{Christian Seis\thanks{Institut f\"ur Angewandte Mathematik, Universit\"at Bonn. Email: seis@iam.uni-bonn.de}}

\maketitle

\begin{abstract}
This review paper is concerned with the stability analysis of the continuity equation in the DiPerna--Lions setting in which the advecting velocity field is Sobolev regular.  Quantitative estimates for the equation were derived only recently \cite{Seis16a}, but optimality was not discussed. In this paper, we revisit the results from \cite{Seis16a}, compare the new estimates with previously known estimates for Lagrangian flows, e.g.\ \cite{CrippaDeLellis08}, and finally demonstrate how those can be applied to produce optimal bounds in applications from physics, engineering or numerics.
\end{abstract}

\section{Introduction}

The linear continuity equation is one of the most elementary partial differential equations. It describes the conservative transport of a quantity by a vector field. We will study this equation in a bounded (Lipschitz) domain\footnote{All the results in this work can be extended to the periodic torus or all of $\R^d$ with suitable modifications.} $\Omega$ in $\R^d$ and denote by $\rho(t,x)\in\R$ and $u(t,x)\in \R^d$ the quantity and the vector field, respectively. For a given initial configuration $\bar \rho(x)\in \R$, the Cauchy problem for the continuity equation reads
\begin{equation}
\label{5}
\left\{
\begin{array}{rll} \partial_t \rho + \div \left(u\rho\right)\! \!\! &= 0&\quad\mbox{in }(0,\infty)\times \Omega,\\
\rho(0,\cdot ) \! \!\!&=\bar \rho &\quad\mbox{in }\Omega.
\end{array}
\right.
\end{equation}
If the vector field is tangential at the boundary of $\Omega$, which we assume from here on, the quantity $\rho$ is (formally) conserved by the flow:
\[
\forall t>0:\quad \int_{\Omega} \rho(t,x)\, dx = \int_{\Omega} \bar \rho(x)\, dx.
\]
Despite its simplicity, the continuity equation plays an important role in fluid dynamics and the theory of conservation laws. In typical applications, $\rho$ represents mass or number density, temperature, energy or phase indicator. In the following, we will frequently refer to $\rho$ as a (possibly negative) mass density, or simply a density. Notice that the vector field has dimensions of length per time, and we will accordingly often refer to $u$ as a velocity field. In fluids applications, $u$ is the velocity of the fluid. 


There is a close link between the partial differential equation (PDE) \eqref{5} and the ordinary differential equation (ODE)
\begin{equation}
\label{1}
\left\{\begin{array}{rl} \partial_t \phi(t,x) &\! \!\!= u(t,\phi(t,x)),\\ \phi(0,x) &\!\!\!=x.\end{array}\right.
\end{equation}
While the PDE represents the Eulerian specification of the flow field, i.e., the description of the dynamics at a fixed location and time, the ODE is the Lagrangian specification, which traces single particles through space and time. The two specifications are in fact equivalent: In the smooth setting, the solution to the continuity equation takes on the nice form
\begin{equation}\label{24}
\rho(t,\phi(t,x))\det \grad \phi(t,x) = \bar \rho(x),\quad\mbox{or simply} \quad\rho(t,\tacka) = (\phi(t,\tacka))_{\#}\bar \rho,
\end{equation}
that is, $\rho$ is the push-forward of $\bar \rho$ by the flow $\phi$,
and a similar formula holds true in the non-smooth setting --- as long as \eqref{5} and \eqref{1} are well-posed. 
This superposition principle is reviewed in \cite[Section 3]{AmbrosioCrippa14}. Notice that for any fixed time $t$, the mapping $\phi(t,\cdot )$ is a diffeomorphism on $\Omega$, whose existence is obtained by the classical Picard--Lindel\"of theorem, and $\det \grad\phi(t,\cdot)$ is the Jacobian determinant, that will be denoted by $J \phi(t,\tacka)$ in the following. The solution $\phi$ of \eqref{1} is called the \emph{flow} of the vector field $u$. 

Out of the smooth setting, well-posedness theory for both the PDE \eqref{5} and the ODE \eqref{1} is more challenging. We focus on the continuity equation from here on and we start with a suitable concept of generalized solutions in the case where $\bar \rho \in L^q(\Omega)$  for some  $q\in[1,\infty]$. We call $\rho$ a \emph{distributional solution} to the continuity equation \eqref{5} in the time interval $[0,T]$, if it conserves the integrability class of the initial datum, $\rho\in L^{\infty}((0,T);L^q(\Omega))$, and satisfies
\[
-\int_0^T \int_{\Omega} (\partial_t \zeta +u\cdot \grad \zeta) \rho\, dxdt = \int_{\Omega} \bar\rho\, \zeta(0,\cdot)\, dx
\]
for any $\zeta\in C^{\infty}_c([0,T)\times \Omega)$. This distributional formulation is reasonable if  $u\in L^1((0,T);L^p(\Omega))$ with $1/p+1/q=1$.

In order to prove existence of distributional solutions, we shall impose a condition on the compressibility of the vector field: If $u$ is  
 weakly compressible,
\begin{equation}
\label{30}
(\div u)^-\in L^1((0,T);L^{\infty}(\Omega)),
\end{equation}
existence is easily obtained by approximation with smooth functions. Here we have used the superscript minus sign to denote the negative part of the divergence.

The questions of uniqueness and continuous dependence on the initial data are more delicate and have been first answered positively by DiPerna and Lions in their ground breaking paper \cite{DiPernaLions89}. Their theory is based on a new solution concept, the theory of renormalized solutions. A \emph{renormalized solution} is a distributional solution $\rho$ with the property that  for any bounded function $\beta\in C^1(\R)$ with bounded derivatives and $\beta(0)=0$ the composition $\beta(\rho)$ satisfies the continuity equation with source
\[
\partial_t\beta(\rho) + \div(u\beta(\rho)) = (\div u)\left(\beta(\rho ) - \rho \beta'(\rho)\right)
\]
in the sense of distributions. In fact, under the additional assumption that  $u\in L^1((0,T);W^{1,p}(\Omega))$, DiPerna and Lions show that distributional solutions \emph{are} renormalized solutions. (The result has been later extended by Ambrosio to vector fields with bounded variation ($BV$) regularity \cite{Ambrosio04} and recently by Crippa, Nobili and the author to the case where velocity gradient is given by a singular integral of an $L^1$ function \cite{CrippaNobiliSeis16}.) The advantage of this solution concept is apparent: By choosing $\beta(s) $ as a suitable approximation of $ |s|^q$, we obtain by integration over $\Omega$ that
\[
\frac{d}{dt} \int_{\Omega} |\rho|^q\, dx = -(q-1) \int_{\Omega} (\div u) |\rho|^q\, dx \le (q-1) \|(\div u)^{-}\|_{L^{\infty}(\Omega)}  \int_{\Omega} |\rho|^q\, dx,
\]
and thus with the help of the Gronwall lemma,
\begin{equation}
\label{10}
\sup_{(0,T)} \|\rho\|_{L^q(\Omega)} \le \exp\left(\int_0^T \|(\div u)^-\|_{L^{\infty}(\Omega)}\, dt\right)^{1-\frac1q} \|\bar\rho\|_{L^q(\Omega)}.
\end{equation}
By the linearity of the continuity equation, this estimate  implies both uniqueness and continuous dependence on the initial data. 

Besides proving well-posedness, DiPerna and Lions study stability under approximations of the vector fields and under diffusive perturbations of the equation. (Notice that this gives two different ways of regularizing the PDE.) While qualitative stability estimates are obtained easily via renormalization, the theory fails to provide quantitative stability estimates that capture the rate of convergence of approximate or perturbative solutions to the original one.  Such estimates were  recently developed in \cite{Seis16a}. 

The aim of the present paper is to revisit the stability estimates from  \cite{Seis16a} and to reformulate them in a new and optimal way. We will mainly focus on two aspects: Our first intention is to compare the new results with earlier achievements in the theory of Lagrangian flows \cite{CrippaDeLellis08} (that, in fact, strongly inspired the estimates in \cite{Seis16a}). Doing so, we hope to convince the reader that the quantities considered in \cite{Seis16a} appear naturally in the context of continuity equations. Secondly, we will  present  applications of the estimates that allow to compute optimal convergence rates in examples of approximate vector fields, zero-diffusivity limits, fluid mixing, and numerical upwind schemes. The last two examples are taken from the studies \cite{Seis13b} and  \cite{SchlichtingSeis16a}. We include these results in order to demonstrate the strength of the estimates from \cite{Seis16a} and to underline the intrinsic connection between the respective works (in particular \cite{Seis13b}) and the latter.
The first example partially extends recent considerations from \cite{DeLellisGwiazdaSwierczewska16}.

We finally like to conclude this introduction by remarking that, as a by-product of the stability estimates, in \cite{Seis16a} a new proof of uniqueness is given for \eqref{5}. This new proof does not rely on the theory of renormalized solutions but is solely based on the distributional formulation of the equation. In a way, the theory in \cite{Seis16a} is the PDE counterpart of the quantitative theory for Lagrangian flows developed by Crippa and De Lellis in \cite{CrippaDeLellis08}. In fact, some of the key estimates were successively transferred from \cite{CrippaDeLellis08} to \cite{Seis16a}.

%

\medskip

{\bf Notation}: In the following, we will use the shorter notation $L^r $ for the Lebesgue space $L^r(\Omega)$,  and similarly $L^r(L^s)$ for $L^r((0,T);L^s(\Omega))$. Further function spaces like $L^1(W^{1,p})$ are defined analogously.

We will omit the domain of integration in the spatial integrals for notational convenience. For instance, we write  $\int\cdot \, dx$ for $ \int_{\Omega}\cdot\, dx$.

We use the sloppy notation $a\lesssim b$ if $a \le C b$ for some constant $C$ that may only depend on the dimension $d$, the domain $\Omega$ or the Sobolev exponent $p$. We write $a\lesssim_{r_1,\dots, r_n} b$ if $C$ depends in addition on the quantities $r_1,\dots,r_n$. Finally, we will sometimes use the notation $a\sim b$ if $a\lesssim b$ and $b\lesssim a$.

\section{Stability estimates for Lagrangian flows}

To motivate our new perspective on the results from \cite{Seis16a}, we start with recalling some facts from the theory of particles moving in a weakly compressible fluid. 
The trajectory of a particle moving with the flow  is given by the solution of the ODE \eqref{1}.
In the classical setting, when the advecting velocity field $u$ is smooth or at least Lipschitz continuous in the spatial variable, existence and uniqueness of a solution is provided by the Picard--Lindel\"of theorem. The Lipschitz regularity also yields simple estimates on the distance of particle trajectories at any time during the evolution.  Indeed, as a consequence of the elementary computation
\[
\left|\frac{d}{dt}  |\phi(t,x) - \phi(t,y)| \right| \le |u(t,\phi(t,x)) -u(t,\phi(t,y))| \le \|\grad u\|_{L^{\infty}}|\phi(t,x) - \phi(t,y)|
\]
and the Gronwall lemma, we easily derive the estimate
\begin{equation}
\label{1a}
\exp\left(- \int_0^t\|\grad u\|_{L^{\infty}}\, dt \right) \le \frac{|\phi(t,x) - \phi(t,y)|}{|x-y|} \le\exp\left(\int_0^t \|\grad u\|_{L^{\infty}}\, dt\right).
\end{equation}
Here we have used Rademacher's identification of Lipschitz functions  with the Sobo\-lev class $W^{1,\infty}$, so that $\|\grad u\|_{L^{\infty}}$ is the Lipschitz constant of $u$.
This estimate illustrates the well-known fact that two particles transported by the flow can neither converge nor diverge faster than exponentially in time.

This classical result can be equivalently rewritten as
\begin{equation}
\label{1b}
- \int_0^t\|\grad u\|_{L^{\infty}}\, dt \le \log \left(\frac{|\phi(t,x) - \phi(t,y)|}{|x-y|}\right) \le \int_0^t \|\grad u\|_{L^{\infty}}\, dt,
\end{equation}
showing that the velocity gradient controls the  logarithmic relative distance of two particles. Here ``relative distance'' refers to the actual distance of particles relative to their initial distance. We will see in the following that it is the latter perspective rather than the classical one \eqref{1a} that allows for a generalization to  the case of flows for less regular vector fields and also to the Eulerian setting.

Notice that both estimates \eqref{1a} and \eqref{1b} contain some information on the regularity of the flow: The flow itself is spatially  Lipschitz, uniformly in time, with Lipschitz constant depending on the gradient of $u$.

Instead of tracing the distance of two different particles in a fluid, we can similarly study the distance of trajectories corresponding to a particle transported by different vector fields: If $\phi$ and $\phi_k$ denote the flows associated via \eqref{1} to the vector fields $u$ and $u_k$, respectively, were $u_k$ may be thought of as a Lipschitz continuous perturbation of $u$, then a computation similar to the one above yields the estimate
\begin{equation}\label{2}
\log\left(\frac{|\phi(t,x) - \phi_k(t,x)|}{\delta} +1\right) \le \int_0^t\|\grad u\|_{L^{\infty}}\, dt  + \frac1{\delta}\int_0^t \|u-u_k\|_{L^{\infty}}\, dt,
\end{equation}
for any $\delta>0$. Therewith, choosing $\delta =\delta_k(t) = \int_0^t \|u-u_k\|_{L^{\infty}}\, dt$, we see that
\begin{equation}
\label{3}
\log\left(\frac{|\phi(t,x) - \phi_k(t,x)|}{\delta_k(t)} +1\right) \le \int_0^t\|\grad u\|_{L^{\infty}}\, dt +1,
\end{equation}
so that, as before, the velocity gradient controls the  logarithmic relative distance of particles moving with two different flows. Observe that $\delta_k(t)$ scales like a length, and can thus be interpreted as the (maximal) distance between the flow fields. Hence, opposed to the situation in \eqref{1b}, we control the distance of particles relative to the distance of vector fields.

Inequality \eqref{3} is an estimate on the rate of convergence of trajectories associated with the vector fields $u$ and $u_k$, if the approximating vector field $u_k$ converges to $u$ in the sense that $\delta_k(t)\to 0$. The statement then shows that the particle trajectories approach each other with a rate of at least $\delta_k(t)$.

Notice that \eqref{2} also implies uniqueness of \eqref{1} when the existence of a solution to the ODE is known. Indeed, if $u=u_k$ is spatially Lip\-schitz and $\phi$ and $\phi_k$ are two solutions to \eqref{1}, the right-hand side of \eqref{2} is bounded independently of $\delta$. Hence, choosing $\delta$ arbitrarily small we see that $\phi $ and $ \phi_k$ must be identical.

Out of the smooth setting, the notion of flows for vector fields has to be appropriately generalized. A common generalization is the notion of regular Lagrangian flows that are well-defined if $u$ is merely Sobolev (or even $BV$) regular in the spatial variable and weakly compressible \cite{DiPernaLions89,Ambrosio04,CrippaDeLellis08}. The latter is expressed by the requirement that
\[
-\infty<   \div u(t,x)\quad \mbox{for a.e.\ }(t,x)\in (0,T)\times \Omega,
\]
cf.\ \eqref{30},
which in turn implies that the Jacobi determinant is bounded below:
\begin{equation}
\label{11}
 J\phi(t,x) = \exp \left(\int_0^t \div u (t,\phi(t,x))\, dt\right) \ge\exp \left(-\int_0^t \|(\div u)^{-} \|_{L^{\infty}}\, dt\right)  = :\Lambda.
\end{equation}
The weak compressibility condition excludes the possibility of infinitely strong sinks in which particles collide in finite time.

Existence, uniqueness and stability of regular Lagrangian flows have been established by DiPerna and Lions in their seminal paper \cite{DiPernaLions89} in the case of vector fields with spatial Sobolev regularity (under the assumption that the divergence is uniformly bounded). This theory has been substantially extended to $BV$ vector fields by Ambrosio \cite{Ambrosio04}. We refer the interested reader to the papers \cite{DeLellis08a,DeLellis08b,AmbrosioCrippa14} for more details and further references, and remark in addition that a comprehensive analysis of the Jacobian is contained in \cite{ColomboCrippaSpirito15}.

Interestingly, DiPerna's and Lions's theory for the ODE \eqref{1} is built on a well-posedness theory for the associated transport  (cf.\ \eqref{23} below) and continuity equations, that is, on the Eulerian (and thus PDE) perspective on particle dynamics. The drawback of the qualitative theory is that no
quantitative estimates can be provided. 
Stability estimates of the type \eqref{3} in the DiPerna--Lions setting were derived later by Crippa and De Lellis \cite{CrippaDeLellis08}, which are of the form
\begin{equation}
\label{4}
\avint \log\left(\frac{|\phi(t,x) - \phi_k(t,x)|}{\delta_k(t)} +1\right)\, dx \lesssim_{\Lambda} \int_0^t\|\grad u\|_{L^{p}}\, dt +1,
\end{equation}
where now
\[
\delta_k(t) = \int_0^t \|u-u_k\|_{L^p}\, dt.
\]
We have decided to work with averaged spatial integrals in all formulas in this paper. With that, $\int_0^t \|\grad u\|_{L^p}\, dt$ is dimensionless and  $\delta_k(t)$ scales like a length, which we interpret, as before, as the distance between the vector fields $u$ and $u_k$.

The papers \cite{DeLellis08a,DeLellis08b,AmbrosioCrippa14} provide reviews of DiPerna's and Lions's theory and of Crippa's and De Lellis's contribution.

Obviously, the above result  confirms that the earlier principle remains valid: Also in the Sobolev case does the velocity gradient provide control over the logarithmic relative distance of particle trajectories. Moreover, the rate of convergence of the trajectory $\phi_k$ to the trajectory $\phi$ is at least of order $\delta_k(t)$ if the latter is tending to zero.

In this weaker setting,  the control of the logarithmic distance ceases to hold uniform in space. Nevertheless, the authors are able to deduce local Lipschitz bounds for the generalized flow. (See also \cite{AmbrosioLecumberryManiglia05} for earlier similar results in this direction.) Moreover, uniqueness can be obtained in a way similar to the one outlined above in the case of Lipschitz vector fields.

The quantitative theory of Crippa and De Lellis fails to cover the full range of vector fields considered earlier by DiPerna and Lions \cite{DiPernaLions89} and Ambrosio \cite{Ambrosio04}. Instead, the authors have to restrict the setting to Sobolev regular vector fields $u\in L^1((0,T);W^{1,p}(\Omega))$ with $p>1$. 
The reason for this is of technical nature: The authors cleverly exploit standard tools from harmonic analysis (maximal functions) whose strong properties just cease to hold if $p=1$.
Stability estimates in the case $p=1$ (and also the $BV$ case) are still open. On the positive side, in \cite{Jabin10}, Jabin manages to extend  estimate \eqref{4} modulo to a factor of order $o(|\log\delta|)$ to the $W^{1,1}$ setting. This estimate is still strong enough to yield uniqueness and stability --- but without rates. A direct proof of uniqueness in the $BV$ setting, that means, without using the uniqueness of the associated partial differential equations as in \cite{Ambrosio04}, was obtained (partially) by Jabin \cite{Jabin10}  and by Hauray and Le Bris \cite{HaurayLeBris11}. A further extension to the case where the velocity gradient is given by a singular integral is treated by Bouchut and Crippa \cite{BouchutCrippa13}.

It remains to remark that stability estimates in the case $p=1$ are closely related
 to a mixing conjecture by Bressan \cite{Bressan03}.  Indeed, Crippa and De Lellis derive the $p>1$ analogue of this conjecture in their paper \cite{CrippaDeLellis08} from an estimate similar to \eqref{4}. See also Subsection  \ref{S:Mixing} (or references \cite{Seis13b,IyerKiselevXu14}) for the corresponding result in the Eulerian setting.

\section{Stability estimates for continuity equations}

 In this section, we will present stability estimates  in the Eulerian framework that are similar to the ODE theory in \cite{CrippaDeLellis08}. That is, instead of tracing single particles in a fluid, we will study the evolution of macroscopic density functions. Our first intention here is to work out analogies to the Lagrangian framework. For this purpose, we study the case of approximate vector fields in Subsection \ref{perturb flow}. Like the estimates in \eqref{3} and \eqref{4}, the result   will be quite general as no relation between the two advecting velocity fields is assumed.  In Subsection \ref{S:diffusion}, we study convergence rates for the zero-diffusivity limit. Subsection \ref{upwind} is devoted to the  convergence order of the numerical upwind scheme. We conclude this paper with an estimate on mixing rates in Subsection \ref{S:Mixing}. We start with the introduction of some notation.

%

\subsection{Kantorovich--Rubinstein distance}\label{KR distance}
In order to transfer the Lagrangian stability estimate \eqref{4} to the Eulerian specification we need some preparations. The quantity that will replace Crippa's and De Lellis's logarithmic trajectory distance is a Kantorovich--Rubinstein distance with logarithmic cost function taken from the theory of optimal transportation and given by
\begin{equation}
\label{10a}
\D_{\delta}(\rho_1,\rho_2) = \inf_{\pi\in \Pi(\rho_1,\rho_2)} \aviint \log\left(\frac{|x-y|}{\delta}+1\right) d\pi(x,y) .
\end{equation}
Functionals of this type were initially introduced to model  minimal costs for transporting mass from one configuration to the other. For two nonnegative distributions $\rho_1$ and $\rho_2$, the set $\Pi(\rho_1,\rho_2)$ consists of all transport plans $\pi$ that realize this transport, i.e.,
\[
\pi[A\times \Omega ] = \int_A \rho_1\, dx,\quad\pi[\Omega\times A] = \int_A\rho_2\, dx,
\]
for any measurable set $A$. The integrand in \eqref{10a} is the so-called cost function that determines the price for the transport between two points\footnote{Concave cost functions are indeed natural in economics applications as they allow to incorporate the \emph{economy of scale} into the mathematical model.}. We refer the interested reader to Villani's monograph \cite{Villani03} for a comprehensive introduction into this topic.

In order to compare this Kantorovich--Rubinstein distance to the trajectory distance considered by Crippa and De Lellis, we notice that in the case where $\rho_1$ and $\rho_2$ can be written as push-forwards of the same configuration, which is, for instance, the case if $\rho_1$ and $\rho_2$ are advected by different flow fields $\phi_1$ and $\phi_2$ while  having the same initial configuration $\bar \rho$ (cf.\ \eqref{24}), then $d\pi = (\phi_1\otimes \phi_2)_{\#} \delta_{x=y}\otimes d\bar \rho$ defines an admissible transport plan in $\Pi(\rho_1,\rho_2) = \Pi((\phi_1)_{\#}\bar \rho,(\phi_2)_{\#}\bar \rho)$. In particular,
\begin{equation}
\label{8}
\D_{\delta}(\rho_1,\rho_2) \le \avint \log\left(\frac{|\phi_1(x)-\phi_2(x)|}{\delta}+1\right) \bar\rho(x)\, dx,
\end{equation}
which means that the Kantorovich--Rubinstein distance $\D_{\delta}(\rho_1,\rho_2)$  is controlled by a weight\-ed variant of Crippa's and De Lellis's logarithmic trajectory distance.

Let us now review some of the properties which make 
 Kantorovich--Rubinstein distances convenient in the study of stability estimates for continuity equations. In fact, what is of importance in our theory is that the quantity $\D_{\delta}(\rho_1,\rho_2)$ 
constitutes a mathematical {\em distance} on the space of configurations of same total mass \cite[Theorem 7.3]{Villani03}, and it metrizes weak convergence \cite[Theorem 7.12]{Villani03}. That is,
\[
\D_{\delta}(\rho_k,\rho)\longrightarrow 0\quad\Longleftrightarrow\quad \rho_k\longrightarrow 0\mbox{ weakly.}
\]
If there exists a sequence of $\delta_k$'s decaying to zero as $k\to \infty$ and such that $\D_{\delta_k}(\rho_k,\rho)$ is uniformly bounded, the latter thus yields that $\rho_k$    converges weakly to $\rho$ with a rate not larger than $\delta_k$.

A crucial insight the stability analysis of \cite{Seis16a} is based on is the dual formulation brought to us in the Kantorovich--Rubinstein theorem
\[
\D_{\delta}(\rho_1,\rho_2) = \sup_{\zeta}\left\{\avint \zeta(\rho_1-\rho_2)\, dx: \: |\zeta(x)-\zeta(y)|\le \log\left(\frac{|x-y|}{\delta}+1\right)\right\},
\]
cf.\ \cite[Theorem 1.14]{Villani03}. One of its immediate consequences is that $\D_{\delta}(\rho_1,\rho_2)$ is a transshipment cost which only sees the difference between $\rho_1$ and $\rho_2$ (``shared mass stays in place''). We can thus write
\[
\D_{\delta}(\rho_1,\rho_2) = \D_{\delta}(\rho_1-\rho_2)\quad\mbox{or}\quad\D_{\delta}(\rho) = \D_{\delta}(\rho^+,\rho^-),
\]
if $\rho^+$ and $\rho^-$ denote, respectively, the positive and the negative part of $\rho$. It follows that Kantorovich--Rubinstein distances can be considered as distances between any two \emph{not necessarily nonnegative} configurations with same average.

In \cite{Seis16a}, the author computed the rate of change of the Kantorovich--Rubinstein distance under the continuity equation with source:
\[
\partial_t \rho + \div(u\rho)   = \div \sigma.
\]
Extending some of the techniques that were developed earlier in \cite{CrippaDeLellis08,BOS,OSS,Seis13b} in the Lagrangian setting, he found that
\begin{equation}
\label{12}
\left|\frac{d}{dt} \D_{\delta}(\rho) \right| \lesssim_{\Lambda,\bar \rho} \|\grad u\|_{L^p}  + \frac1{\delta} \|\sigma\|_{L^1},
\end{equation}
cf.\ \cite[Proposition 1]{Seis16a}, if $\rho$ has zero mean.

\subsection{Approximating the vector field}\label{perturb flow}
We now consider the situation from \eqref{4}. 
We thus let $u$ and $u_k$ be two vector fields in $L^1(W^{1,p})$ satisfying the compressibility condition \eqref{11}, and we denote by $\rho$ and $\rho_k$ the corresponding solutions to the continuity equation \eqref{5}, starting with the same initial datum $\bar \rho$ in $L^q$ with $1/p + 1/q=1$. As an immediate consequence of \eqref{12} and \eqref{10}, we obtain our first result.

\begin{theorem}[\cite{Seis16a}]\label{T1}
If  $\delta_k(t)$ denotes the distance between the vector fields $u$ and $u_k$ given  by
\[
\delta_k(t) = \int_0^t \|u-u_k\|_{L^p}\, dt,
\]
then it holds
\[
\D_{\delta_k}(\rho,\rho_k) \lesssim_{\Lambda,\bar \rho} \int_0^t \|\grad u\|_{L^p}\, dt +1.
\]
\end{theorem}

Notice that there is similarity to the control principle we found earlier in the ODE case: The velocity gradient controls the logarithmic relative distance of two configurations.

With regard to the fact that Kantorovich--Rubinstein distances metrize weak convergence, in the situation where $\delta_k\to0$, the statement in the theorem now shows that
\[
\rho_k\longrightarrow \rho \mbox{ weakly with rate not larger than }   \delta_k.
\]
This estimate is sharp as can be seen by the following example suggested by De Lellis, Gwiazda and \'Swierczewska-Gwiazda \cite{DeLellisGwiazdaSwierczewska16}. 

\begin{example}
Consider the oscillating vector field  $u_k(x) =  \sin(2\pi kx)/2\pi k$ on the interval $\Omega = [0,1]$. Solving the continuity equation with the initial datum $\bar \rho=1$ yields the oscillating solution 
\[
\rho_k(t,x) = \frac{1+ \tan^2(\pi k x)}{e^t +e^{-t} \tan^2(\pi k x)},
\]
cf.\ Figure \ref{fig1}. Because $u_k$ converges strongly to zero as $k \to \infty$, it is clear that the limiting problem is stationary, i.e., $\rho\equiv 1$.
\begin{figure}[t]
\begin{center}
\includegraphics[width=.6\textwidth]{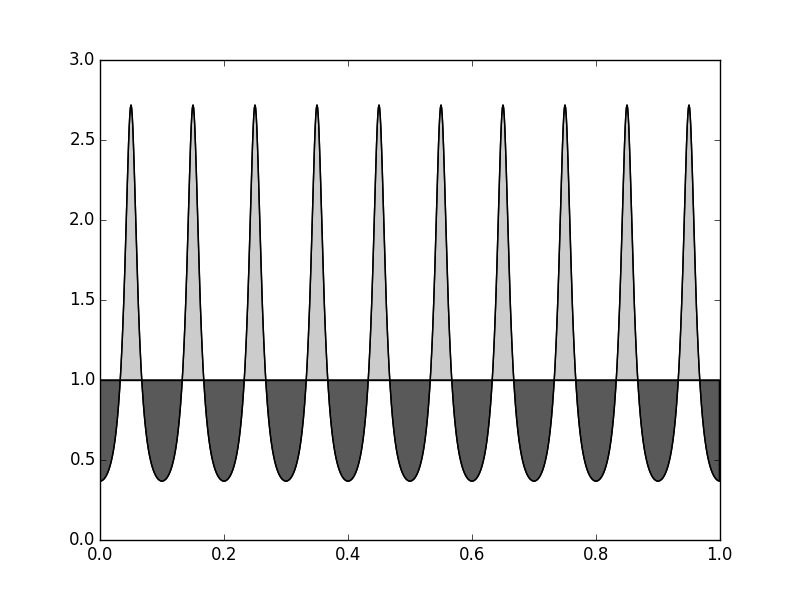}
\end{center}
\caption{The figure shows the oscillating density $\rho_{10}$ at time $t=1$. The corresponding Kantorovich--Rubinstein distance $\D_{\delta_{10}}(\rho_{10},\rho)$ measures the transport between $(\rho_{10}-\rho)^+$ (lightgrey region) and $(\rho_{10}-\rho)^-$ (darkgrey region).}\label{fig1}
\end{figure}

In view of the oscillatory behavior of $\rho_k$, the convergence to $\rho\equiv 1$ is merely weakly,
\[
\|\rho-\rho_k\|_{L^1(L^1)} \noto 0.
\]
In order to quantify the rate of weak convergence, we notice that
\[
\delta_k(t)  = \int_0^t \|u-u_k\|_{L^p}\, dt = \frac{t}{2\pi k} \left(\int_0^1 |\sin(2\pi k x)|^p\, dx\right)^{1/p} \sim \frac{t}{k}
\]
because $u=0$. By the periodicity and the symmetry of the problem, we furthermore compute 
\[
\D_{\delta_k}(\rho,\rho_k)  = k \D_{\delta_k}({\rho}_{|_{\left[0,1/k\right]}},  {\rho_k}_{|_{\left[0,1/k\right]}}) = \D_{\delta_1}(\rho,\rho_1) \sim_t 1,
\]
where we have rescaled length in the last identity.
\end{example}

This example shows that Theorem \ref{T1} is sharp in two respects: Firstly, strong convergence of $\rho_k$ to $\rho$ does in general not hold. Secondly, the result captures the correct rate of convergence.
%
%

Even though the result in Theorem \ref{T1} was already present in \cite{Seis16a}, the author deduced a weaker stability estimate in order to replace the unwieldy Kantorovich--Rubinstein distance by a standard negative Sobolev norm. In fact, the author proved that
\[
\|\rho-\rho_k\|_{W^{-1,1}} \lesssim_{\Lambda,\bar \rho,u} \frac1{|\log\delta_k(t)|},
\]
see \cite[Theorem 2]{Seis16a}. The new formulation in Theorem \ref{T1} has the advantage that it is sharp and naturally extends the analogous estimates in the Lagrangian setting \eqref{4}.

\subsection{The zero-diffusivity limit}\label{S:diffusion}

In this subsection, we expand the model \eqref{5} by a second parallel transport mechanism besides advection: diffusion. Advection-diffusion models are ubiquitous in thermodynamics, fluid dynamics and engineering, for instance in the context of thermal convection \cite{Siggia94}, spinodal decomposition \cite{Siggia79}, or mixing \cite{Thiffeault12}. While convection enhances the efficient transport of particles or fluid parcels over large distances and tends to create sharp gradients in the density (or temperature) distribution, diffusion compensates  density (or temperature) differences locally.

We thus consider in the following the Cauchy problem for the advection-diffusion equation
\begin{equation}
\label{6}
\left\{
\begin{array}{rll}
\partial_t \rho_{\kappa} + \div\left(u\rho_{\kappa}\right) \! \!\!& = \kappa\laplace \rho_{\kappa}&\quad\mbox{in }(0,\infty)\times \Omega,\\
\rho_{\kappa}(0,\cdot ) \! \!\!&=\bar \rho &\quad\mbox{in }\Omega,
\end{array}
\right.
\end{equation}
where $\kappa$ is the (positive) diffusivity constant. Equipping the equation with the no-flux condition $\grad\rho_{\kappa}\cdot \nu = 0$ on the boundary of $\Omega$ implies that the evolution is still mass is conserving, 
\begin{equation}
\label{9}
\int\rho_{\kappa}(t,x)\, dx = \int\bar \rho\, dx.
\end{equation}
We assume furthermore that  $\bar \rho$ is nonnegative, and so $\rho_{\kappa}$ is as a consequence of the maximum principle for  \eqref{6}. Notice that \eqref{10} remains valid for $\rho_{\kappa}$, which can be easily seen by testing \eqref{6} with $\rho_{\kappa}^{q-1}$. 

We are interested in the  vanishing diffusivity limit $\kappa\to0$. In order to quantify the rate of convergence of solutions of \eqref{6} towards solutions of the  purely advective model \eqref{5}, we will make use of a  standard decay estimate from the relaxation theory for the diffusion (or heat) equation. A common way to identify the equilibration rate in the diffusive model is by studying the decay behavior of the entropy
\[
H(\rho) = \avint \rho\log \rho\, dx.
\]
We compute the rate of change of entropy under the evolution \eqref{6} using multiple  integrations by parts,
\begin{eqnarray*}
\frac{d}{dt} H(\rho_{\kappa}) &=&  \kappa \avint \laplace \rho_{\kappa} \log\rho_{\kappa}\, dx - \avint \div\left(u\rho_{\kappa}\right)\log\rho_{\kappa}\, dx\\
& =& - \avint \frac{|\grad\rho_{\kappa}|^2}{\rho_{\kappa}}\, dx - \avint (\div u ) \rho_{\kappa}\, dx,
\end{eqnarray*}
where in the first equality we have used the fact that the evolution is mass conserving \eqref{9}.
Moreover, since $\rho_{\kappa}$ is a nonnegative function, integration in time yields
\[
\sup_t H(\rho_{\kappa}) + \kappa \int_0^t \avint \frac{|\grad\rho_{\kappa}|^2}{\rho_{\kappa}}\, dxdt \le H(\bar\rho) + \left(\int_0^t \|(\div u)^{-} \|_{L^{\infty}}\, dt \right)\|\rho_{\kappa}\|_{L^{\infty}(L^1)} .
\]
Then, if the initial density has finite entropy, by H\"older's inequality and mass conservation \eqref{9},
\[
\int_0^t \avint |\grad\rho_{\kappa}|\, dxdt  \le \int_0^t\left( \avint \rho_{\kappa}\, dx  \avint \frac{|\grad\rho_{\kappa}|^2}{\rho_{\kappa}}\, dx\right)^{1/2}dt \lesssim_{\bar \rho,\Lambda} \sqrt{\frac{t}{\kappa}}.
\]
Using the theory developed in \cite{Seis16a}, cf.\ \eqref{12}, and the a priori estimate \eqref{10}, it can now be shown that
\[
\frac{d}{dt}\D_{\delta}(\rho,\rho_{\kappa}) \lesssim_{\Lambda,\bar \rho} \|\grad u\|_{L^p} + \frac{\kappa}{\delta}\|\grad\rho_{\kappa}\|_{L^1}.
\]
Integration in time and a combination of the previous two estimates then yields the following result: 

\begin{theorem}\label{T2}
Let $\delta_{\kappa}(t) $ be the diffusion distance per time $t$, i.e., $\delta_{\kappa}(t)= \sqrt{t\kappa}$. Then
\[ 
\D_{\delta_{\kappa}(t)} (\rho,\rho_{\kappa}) \lesssim \int_0^t \|\grad u\|_{L^p} + 1.
\]
\end{theorem}

In other words, the diffusive approximation  converges weakly to the unique solution of the continuity equation with a rate not larger than $\sqrt{t\kappa}$. The latter equals approximately the distance a particle can travel by diffusion in time $t$.
Also in this case, a qualitative convergence result was previously established in the work of DiPerna and Lions \cite{DiPernaLions89}. To the best of our knowledge, it is for the first time that a convergence rate for the zero-diffusivity limit is obtained.

\subsection{Convergence rates for the upwind scheme}\label{upwind}


The upwind scheme is a numerical scheme for approximating solutions to the continuity equation. The scheme is a finite volume scheme, which means, that the domain is decomposed into control volumes (or cells) of small diameter and the evolving density is approximated by averages over each control volume.

To be more specific, we consider a domain $\Omega$ that can be written as a finite union of rectangular boxes. We decompose $\Omega$   into a family of rectangular cells with disjoint interiors, $\Omega = \cup_{K\in\cT} K$, where $\cT$ is the tessellation and $K$ is a translation of the cube $ [0, h_1]\times \dots\times [0,h_d]$.  The size $h$ of the tessellation is the maximal edge length, that is,
\[
h  =  \max_{i=1}^d h_i.
\]
We suppose that the tessellation is regular in the sense that  $h_i\sim h$ for all $i$. For two neighboring cells $K\sim L$, we denote by $K\edge L$ the joint boundary.  The normal vector on $K\edge L$ pointing from $K$ to $L$ is denoted by $\nu_{KL}$.

We choose a fixed time step size $\delta t$  so that the $n$-th time step reads $t^n = n\delta t$. To guarantee the stability of the explicit scheme, 
we impose the following Courant--Friedrichs--Lewy (CFL) condition on the time step size:
\[
\forall n: \quad \int_{t^n}^{t^{n+1}} \|u\|_{L^{\infty}}\, dt \le h.
\]
We are thus assuming in this subsection that $u\in L^1(L^{\infty})$.

To approximate the transport term, we consider the net flow  from $K$ to $L$ defined by
\[
u_{KL}^n  = \avint_{t^n}^{t^{n+1}} \avint_{K\edge L} u\cdot \nu_{KL}\, d\Ha^{d-1} dt.
\]
We remark that these quantities are well-defined thanks to the trace estimate for Sobolev functions. Furthermore, the initial configuration of the scheme is the volume average over each $K\in\cT$, i.e.,
\[
\rho_K^0 = \avint_K \bar \rho\, dx.
\]
We are now in the position to define the explicit upwind finite volume scheme for the continuity equation \eqref{5},
\begin{equation}
\label{7}
\rho_K^{n+1} = \rho_K^n + \frac{dt}{h}\sum_{L\sim K}\left( u_{LK}^{n+}\rho_L^n  - u_{KL}^{n+}\rho_K^n\right),
\end{equation}
where $u_{KL}^{n+} = (u_{KL}^n)^+$. The approximate solution is given by
\[
\rho_h(t,x) = \rho_K^n\quad\mbox{if }(t,x)\in[t^n,t^{n+1})\times K,\, K\in\cT.
\]
See \cite{EymardGallouetHerbin00,LeVeque02} for properties and reference.
In the DiPerna--Lions setting, convergence of the scheme, i.e., $\rho_h\to \rho$ as $h\to 0$, 
 is proved in \cite{Walkington05,Boyer12,SchlichtingSeis16a}.

Even though the numerical scheme is formally first order,  one observes a break down in the convergence rate to order $1/2$ in the case of non-smooth initial data. In the DiPerna--Lions setting considered here, $\sqrt{h}$-rates were numerically observed in \cite{Boyer12} and \cite{SchlichtingSeis16a}.
The reason for this lack of convergence is the occurrence of numerical diffusion that smooths out sharp interfaces. Such irregularities, however, are simply transported in the continuous model. In a certain sense, approximate solutions show a behavior similar to those of the advection-diffusion equation \eqref{6}, where $\kappa\sim h$. It is this similarity that determines the $\sqrt{h}$-rate of convergence, cf.\ Theorem \ref{T2} above. The effect of numerical diffusion is illustrated in Section 2.4 of \cite{SchlichtingSeis16a}. 

In the case of regular (i.e., at least spatially Lipschitz continuous) vector fields, this break down in the order of convergence is long known. First rigorous results on optimal convergence rates date back to the 1970s, see e.g., \cite{Kuznecov76,Peterson91,VilaVilledieu03,Despres04,
MerletVovelle07,Merlet07,CockburnDongGuzmanQian10,
DelarueLagoutiere11,DelarueLagoutiereVauchelet16}. To the best of our knowledge, the only available result in the DiPerna--Lions setting is very recent: In \cite{SchlichtingSeis16a}, the author establishes jointly with Schlichting an upper bound on the rate of weak convergence that captures the optimal order.

\begin{theorem}[\cite{SchlichtingSeis16a}]
\label{T3}
Let $\delta_h(t)$ be the numerical diffusion distance per time $t$, i.e., $\delta_h(t) = \sqrt{h \int_0^t \|u\|_{L^{\infty}}\, dt}$. Then
\[
\D_{\delta_h}(\rho,\rho_h)\lesssim_{\Lambda,\bar \rho} \int_0^t\|\grad u\|_{L^p}\, dt +1.
\]
\end{theorem}

The work \cite{SchlichtingSeis16a} builds up not only on the quantitative theory from \cite{Seis16a}. A crucial ingredient is a probabilistic interpretation of the upwind scheme suggested by Delarue, Lagouti\`ere and Vauchelet \cite{DelarueLagoutiere11,DelarueLagoutiereVauchelet16}. In fact, in \cite{SchlichtingSeis16a}, we interpret \eqref{7} as a Markov chain, which comes as a time-discretized version of the stochastic differential equation
\[
d\psi_t  = u(t,\psi_t)dt + \sqrt{2h }\, dW_t,
\]
with a noise term depending on the details of the mesh. In a certain sense, the above equation is the Lagrangian analogue of the advection-diffusion equation \eqref{6}. It turns out that the noise term determines the $\sqrt{h}$-rate of convergence.

\subsection{Mixing by stirring}\label{S:Mixing}

In the past years, mixing by stirring attracted much interest in both the applied mathematics and the engineering communities. Mixing refers to the homogenization process of an inhomogeneous substance being stirred by an agent. One of the major goals is the quantification of mixing rates and the design of mixing strategies. In order to optimize mixing strategies, absolute lower bounds on the mixing rate are indispensable. In this subsection, we present a lower bound on mixing by stirring of incompressible viscous fluids that was obtained earlier by the author in \cite{Seis13b}. A  nice review on the mathematical side of mixing was written by Thiffeault \cite{Thiffeault12}.

A natural constraint in the experimental mixing set-up is the amount of mechanical work the engineer is willing to spend in order to overcome viscous friction  to maintain stirring. Mathematically, this amounts to limiting the budget of the viscous dissipation rate (or enstropy) given by $\|\grad u\|_{L^2}$. In the following, we will slightly generalize this constraint by assuming that $u\in L^1(W^{1,p})$ for some $p>1$ as in the previous part of this paper.

While our intuition is strong about whether a substance is well mixed or not, the choice of a measure that quantifies the degree of mixedness depends on the mathematical communities. Among fluid dynamists, homogeneous negative Sobolev norms are favored, in particular the $\dot H^{-1/2}$ norm \cite{MathewMezicPetzold05,MathewMezicGrivopoulos07} and the $\dot H^{-1}$ norm \cite{DoeringThiffeault06,ShawThiffeaultDoering07,LinThiffeaultDoering11}. These norms measure oscillations: the larger the length scales, the larger the negative Sobolev norms. 

In \cite{Seis13b}, the author introduces a new mixing measure besides the $\dot H^{-1}$ norm: a variant of the Kantorovich--Rubinstein distance introduced earlier in this paper. We accordingly consider
\[
M(\rho) = \inf_{\pi\in\Pi(\rho^+,\rho^-)} \exp\left(\aviint \log|x-y|\, d\pi(x,y)\right).
\]
Notice that $M(\rho) = \lim_{\delta\to 0}\exp\left(\D_{\delta}(\rho)  + (\log \delta)\|\rho\|_{L^1}\right)$. In the case of a two-phase mixture, modeled by $\rho \in\{\pm1\}$, this distance formally scales as a length, so $M(\rho)$ agrees with the average size of the unmixed regions.

The mixing process can be modeled by the continuity equation \eqref{5}, which turns into the transport equation
\begin{equation}
\label{23}
\partial_t\rho + u\cdot \grad \rho=0
\end{equation}
under the assumption that the fluid is incompressible $\div u=0$, which we shall assume for convenience. For simplicity, we restrict our attention to two-phase mixtures with equal volume fraction, so that 
\[
|\left\{x\in\Omega: \rho(t,x)=1\right\}| = |\left\{x\in\Omega:\: \rho(t,x)=-1\right\}|
\]
for any $t>0$, or, equivalently, $\int \rho\, dx=0$.

In \cite{Seis13b}, the author derives a lower bound on mixing rates in   incompressible viscous fluids, building up on an estimate similar to \eqref{12}.

\begin{theorem}[\cite{Seis13b}]
\label{T4}
For any $T\ge 0$, it holds that
\begin{equation}\label{22}
M(\rho(T,\tacka)) \ge M(\bar \rho) \exp\left(-\frac1C\int_0^T \|\grad u\|_{L^p}\, dt\right),
\end{equation}
where $C$ is a constant depending only on $p$ and $d$.
\end{theorem}

This estimate shows the impossibility of perfect mixing, i.e., $\rho\to0$ weakly, in finite time. A similar statement has been obtained earlier by Crippa and De Lellis \cite{CrippaDeLellis08} for a certain geometric mixing measure suggested by Bressan \cite{Bressan03}. In fact, Bressan conjectures the $p=1$ analogue of Crippa's and De Lellis's estimate.

It is not difficult to deduce from \eqref{22} a lower bound on the decay rate of the $\dot H^{-1}$ norm. Indeed, in \cite{Seis13b}, it is moreover proved that
\begin{equation}
\label{20}
\frac1{|\grad \rho|_{BV}} \lesssim M(\rho) \le | \grad^{-1}\rho|_{L^2},
\end{equation}
where $|\grad \rho|_{BV}$ and $|\grad^{-1}\rho|_{L^2}$ denote, respectively, the homogeneous part of the $BV$ norm, and thus $|\grad \rho|_{BV} = 2 |\partial\{\rho=1\}|/|\Omega|$, and the homogeneous part of the $H^{-1}$ norm. The first inequality in \eqref{20} is an interpolation inequality, whereas the second one follows immediately via Jensen's inequality and the Kantorovich--Rubinstein theorem \cite[Theorem 1.14]{Villani03}. Plugging \eqref{20} into Theorem \ref{T4} yields
\begin{equation}
\label{21}
|(\grad^{-1} \rho)(T,\tacka) |_{L^2}\gtrsim \frac1{|\grad\rho|_{BV}} \exp\left(-\frac1C\int_0^T \|\grad u\|_{L^p}\, dt\right).
\end{equation}
A similar decay estimate for the $\dot H^{-1}$ norm has been obtained simultaneously by Iyer, Kiselev and Xu \cite{IyerKiselevXu14} by using the geometric results from \cite{CrippaDeLellis08}.

Estimates \eqref{22} ad \eqref{21} are sharp. This was proved  by Yao and Zlato\v{s} \cite{YaoZlatos14} and independently by Alberti, Crippa and Mazzucato \cite{AlbertiCrippaMazzucato14}. In fact, in both works, the authors construct explicit mixing flows that saturate the lower bounds from \cite{Seis13b} and \cite{IyerKiselevXu14}. Numerical evidence for the optimality of this mixing rate was given earlier in \cite{LinThiffeaultDoering11}.

There is a close relation between Theorem \ref{T4} and the lower bound in \eqref{1a}. (In fact, also the upper bound
\[
M(\rho(T,\tacka)) \le M(\bar \rho) \exp\left(\frac1C\int_0^T \|\grad u\|_{L^p}\, dt\right)
\]
is valid.)  Estimate  \eqref{22} can be seen as the Eulerian (and Sobolev) analogue of \eqref{1a}, in the sense that in Theorem \ref{T4}, we compute the distance between the configuration described by the mixing process and the stationary fully mixed state $\rho=0$. While \eqref{1a} shows that trajectories cannot converge faster than exponentially in time, the Eulerian analogue shows that different density configurations cannot converge faster than exponentially in time. This observation also underlines the link between mixing and the question of uniqueness for  the partial differential equation \eqref{5} (or \eqref{23}): A system is perfectly mixing in finite time precisely if solutions to \eqref{5} are in general not unique. Notice that in the case of finite time mixing, upon reversing time, one has nontrivial solutions to \eqref{5} with zero initial datum. An explicit construction of such an unmixing solution is due to Depauw \cite{Depauw03}.

It remains to remark that upper bounds on the rates of unmixing (or coarsening) in viscous fluids were obtained in \cite{BOS,OSS}. The analysis in these papers combines \eqref{22} and the lower bound of \eqref{20} with the Kohn--Otto upper bound method \cite{KohnOtto02}.

\section*{Acknowledgement}
The author thanks Andr\'e Schlichting for fruitful discussions and for suggesting the entropy approach in Subsection \ref{S:diffusion}.
\bibliography{coarsening}
\bibliographystyle{abbrv}

\end{document}